\documentclass[12pt]{amsart}
\usepackage{amssymb, amscd, latexsym, color, amscd, float}
\usepackage{url}
\usepackage[final]{graphicx}

\usepackage{ulem,ifthen}
\newtheorem{thm}{Theorem}[section]

\newtheorem{lem}[thm]{Lemma}
\newtheorem{cor}[thm]{Corollary}

\newtheorem{defn}[thm]{Definition}

\vspace*{-15mm}

\begin{document}

\vspace*{2cm}


\title
[]
{\normalsize A periodicity theorem\\
for acylindrically hyperbolic groups}

\bigskip

\author{Oleg Bogopolski}
\address{{Sobolev Institute of Mathematics of Siberian Branch of Russian Academy
of Sciences, Novosibirsk, Russia}\newline {and D\"{u}sseldorf University, Germany}}
\email{Oleg$\_$Bogopolski@yahoo.com}

\begin{abstract}
We generalize a well known periodicity lemma
from the case of free groups to the case of acylindrically hyperbolic groups.
This generalization will be used later to describe solutions of certain equations in acylindrically hyperbolic groups and to characterize verbally closed finitely generated acylindrically hyperbolic subgroups of
finitely presented groups.

\end{abstract}

\maketitle

\bigskip

\section{Introduction}

\subsection{The main result}
Let $A$ be an alphabet and $u,v$ be two words in the free monoid $A^{\ast}$.
The word $v$ is called a {\it period} of $u$ if there exists $n\geqslant 1$
such that $v^n$ is an initial subword of $u$ and $u$ is an initial subword of $v^{n+1}$.
The length of $u$ in the alphabet $A$ is denoted by $|u|$.

The following periodicity lemma is a kind of folklore that was known to researchers
working in such areas as combinatoric of words and Burnside problems
in 1960's or earlier.

\begin{lem}\label{Theorem_FW}  {\rm } {\rm (Periodicity lemma)}
Let $z$ be a word in a free monoid $A^{\ast}$.\break
If $z$ has periods $u$ and~$v$ such that $|z|\geqslant |u|+|v|$, then $z$ has a period $w$ such that $u=w^s$ and $v=w^t$
for some $s,t\in \mathbb{Z}$.
\end{lem}

Let $F_r$ be the free group of rank $r$ and let $F_r^n$ be its subgroup generated by all $n$-th powers.
The group ${\bf B}_r(n)=F_r/F_r^n$ is called the free Burnside group of rank $r$ and exponent $n$.
In 1968, Novikov and  Adian solved the bounded Burnside problem by proving that ${\bf B}_r(n)$ is infinite for $r\geqslant 2$ and odd $n\geqslant 4381$
(see~\cite{NA_1} and further development in~\cite{A,O,O_1,I,L,DG}).
In the proof they use the notion of periods (of ranks $n\in \mathbb{N}$) and the existence of $p$-aperiodic infinite words, i.e. words which do not contain
a finite subword of the form $u^p$. Explicitly, Lemma~\ref{Theorem_FW} can be found in the book of Adian~\cite{A}
(see statement 2.3 in Chapter I there).
In its optimal form, which we do not present here, this lemma was proved by Fine and Wilf in~\cite{FW}.
The usefulness of this lemma in the study of equations in (semi)groups is pointed out in~\cite[Section~5.4]{BP}.


\newpage

We reformulate Lemma~\ref{Theorem_FW} for free groups in the form which enables to see an analogy
with our Theorem~\ref{acylindric}.


\begin{lem}\label{reformulation} {\rm(reformulation of Lemma~\ref{Theorem_FW})}
Let $a,b$ be two cyclically reduced words in the free group $F$ with basis $X$.
If the bi-infinite words $L(a)=\dots aaa\dots $ and $L(b)=\dots bbb\dots$ have a common subword of length $|a|+|b|$,
then some cyclic permutations of $a$ and $b$ are powers of some word $c$.
\end{lem}


\vspace*{-32mm}
\hspace*{-2mm}
\includegraphics[scale=0.7]{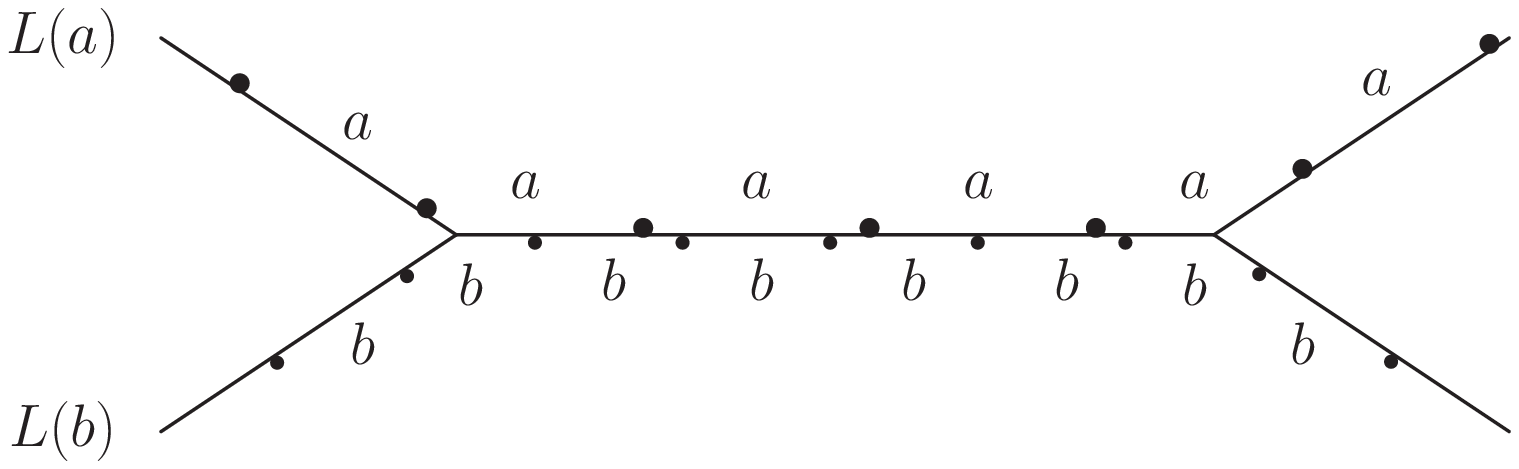}

\vspace*{-14cm}

\begin{center}
Fig. 1. Illustration to Lemma~\ref{reformulation}.
\end{center}


The main result of this paper, Theorem~\ref{acylindric}, generalizes this theorem to
the case of acylindrically hyperbolic groups. 

The class of acylindrically hyperbolic groups
was defined by Osin in~\cite{Osin_2} (see Definition~\ref{definition_acylindrically_hyperbolic} below), see papers
\cite{Osin_2,Osin_3,DOG} for historical aspects, properties and examples.
This class includes non-(virtually cyclic) groups that are hyperbolic
relative to proper subgroups, groups of deficiency at least~2, all but finitely many mapping class groups,
 $Out(F_n)$ for $n\geqslant 2$,
non-(virtually cyclic) groups acting properly on proper CAT(0)-spaces and containing rank-one elements,
non-cyclic directly indecomposable right-angled Artin groups;
see more examples in the survey of Osin~\cite{Osin_3}. 

To formulate the main result, we need some notions.

\begin{defn}\label{L}
{\rm
Let $Y$ be a generating set of a group $G$.
Given two elements $x,a\in G$, we consider the bi-infinite path $L(x,a)$ in the Cayley graph $\Gamma(G,Y)$
obtained by connecting consequent points $\dots, xa^{-1},x,xa,\dots $ by geodesic segments so that,
for all $n\in \mathbb{Z}$, the segment $p_n$ connecting $xa^n$ and $xa^{n+1}$ has the same label as the segment $p_0$ connecting $x$ and $xa$.
The paths $p_n$ are called {\it $a$-periods} of $L(x,a)$. For a subpath $p\subset L(x,a)$  and a number $k\in \mathbb{N}$, we say that the path $p$ {\it contains $k$ $a$-periods} if there exists $n\in \mathbb{Z}$ such that $p_np_{n+1}\dots p_{n+k-1}$ is a subpath of $p$.
The vertices $xa^n$, $n\in \mathbb{Z}$, are called the {\it phase} vertices of $L(x,a)$.
}
\end{defn}



Let $G$ be a group and $X$ a generating set of $G$. Suppose that the Cayley graph $\Gamma(G,X)$
is hyperbolic and that $G$ acts acylindrically on $\Gamma(G,X)$. In~\cite{Bowditch}
Bowditch proved that the infimum of stable norms (see Section 2) of all loxodromic
elements of $G$ with respect to $X$ is a positive number.
We denote this number by ${\bf inj}(G,X)$ and call it the {\it injectivity radius} of $G$ with respect to $X$.


\begin{thm}\label{acylindric}
Let $G$ be a group and $X$ a generating set of $G$. Suppose that the Cayley graph $\Gamma(G,X)$
is hyperbolic and that $G$ acts acylindrically on $\Gamma(G,X)$.
Then there exists a constant $C>0$ such that the following holds.

Let $a,b\in G$ be two loxodromic elements which are shortest in their conjugacy classes
and such that $|a|_X\geqslant |b|_X$.
Let $x,y\in G$ be arbitrary elements and $r$ an arbitrary non-negative integer. We set $f(r)=\frac{2r}{{\bf inj}(G,X)}+C$.

If a subpath $p\subset L(x,a)$ contains at least $f(r)$ $a$-periods and lies in the $r$-neighborhood of $L(y,b)$,
then there exist $s,t\neq 0$ such that $(x^{-1}y)b^s(y^{-1}x)=a^t$.
In particular, $a$ and $b$ are commensurable.
\end{thm}

\medskip

\vspace*{-30mm}
\hspace*{-5mm}
\includegraphics[scale=0.7]{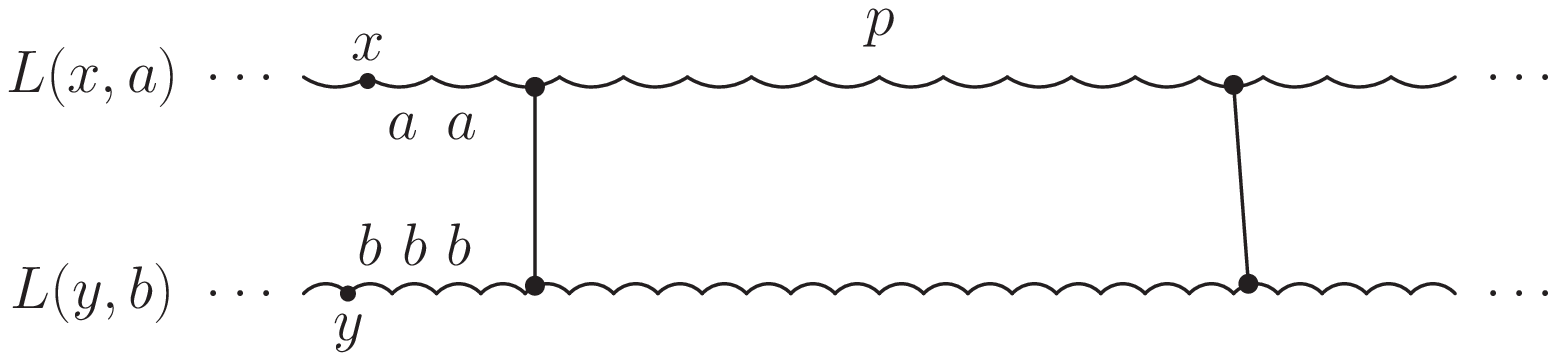}

\vspace*{-14.5cm}

\begin{center}
Fig. 2. Illustration to Theorem~\ref{acylindric}.
\end{center}

\medskip

Observe that if $G$ is a free group, then we can take $f(0)=2$ as follows from Lemma~\ref{reformulation}.

\medskip

\noindent
{\bf Remarks.}
The main issue of Theorem~\ref{acylindric} is that the function $f$ is linear and does not depend on $|a|_X$ and $|b|_X$.
Another point is that $X$ is allowed to be infinite.
Analogous theorem can be also formulated and proved for groups that act acylindrically on hyperbolic spaces 
(not necessarily on their Cayley graphs).

A version of Theorem~\ref{acylindric} (for $r$ smaller than some constant) is contained in the paper of Coulon~\cite[Proposition~3.44]{Coulon}.
Periods were also used 
by Delzant and Steenbock in~\cite{DS} for studying the product set growth in groups and hyperbolic geometry.


\subsection{Applications of the main result.} In~\cite{Bog_1}, we use Theorem~\ref{acylindric}
to describe solutions of certain equations in acylindrically hyperbolic groups and to characterize
finitely generated verbally closed acylindrically hyperbolic subgroups of finitely presented groups.
Below we formulate some of these results. The last one, Corollary~\ref{1_corollary_hyp},
solves a problem of Myasnikov and Roman'kov
(see Problem~5.2 in~\cite{MR}).

\medskip

We recall some definitions; see more details in Section 2.
Suppose that $G$ is a group acting acylindrically on a hyperbolic space.
Then any loxodromic element $g\in G$ is contained in a unique maximal virtually cyclic subgroup $E_G(g)$
of~$G$ \cite[Lemma 6.5]{DOG}.
This subgroup is called the {\it elementary subgroup associated with} $g$.
An element $g\in G$ is called {\it special} if it is loxodromic and $E_G(g)=\langle g\rangle$.
Two elements $a,b\in G$ of infinite order are called {\it commensurable}
if there exist $g\in G$ and $s,t\in \mathbb{Z}\setminus \{0\}$ such that $a^s=g^{-1}b^tg$;
see other definitions in Section~2.


\begin{cor}\label{eq_acylindric} {\rm{(\cite[Corollary 2.8]{Bog_1})}}
Let $G$ be an acylindrically hyperbolic group with respect to a generating set $S$.
Suppose that $a,b\in G$ are two non-commensurable special elements (with respect to $S$).
Then there exists a number $\ell=\ell(a,b)\in \mathbb{N}$ such that for all $n,m\in \mathcal{\ell}\mathbb{N}$, $n\neq m$, the equation $x^ny^m=a^nb^m$ is perfect, i.e. any solution of this equation in
$G$ is conjugate to $(a,b)$ by a power of $a^nb^m$.

\end{cor}

\begin{defn}
{\rm (see~\cite[Definition 1.1]{MR})
A subgroup $H$ of a group $G$ is called {\it verbally closed} in $G$ if
for any word $W$ in variables $x_1,\dots,x_n$, $n\in \mathbb{N}$, and any element $h\in H$ the following holds: if the equation $W(x_1,\dots,x_n)=h$
has a solution in $G$, then it has a solution in $H$.
}
\end{defn}

\begin{thm}\label{1_prop 4.1} {\rm{(\cite[Theorem 2.2]{Bog_1})}}
Suppose that $G$ is a finitely presented group and $H$ is a finitely generated acylindrically hyperbolic subgroup of $G$ such that $H$ does not normalise a nontrivial finite subgroup of $H$.
Then $H$ is verbally closed in $G$ if and only if $H$ is a retract of $G$.
\end{thm}

Recall that a group $H$ is called {\it equationally noetherian} if
every system of equations with constants from $H$ and a finite number of variables
is equivalent to a finite subsystem, see~\cite{BMR}. We say that a group $G$ is {\it finitely generated over a subgroup} $H$ if
there exists a finite subset $A\subset G$ such that $G=\langle A,H\rangle$.


\begin{thm}\label{1_prop 4.2} {\rm{(\cite[Theorem 2.4]{Bog_1})}}
Let $G$ be a group that is finitely generated over a subgroup $H$.
Suppose that $H$ is equationally noetherian, acylindrically hyperbolic, and without a nontrivial finite normal subgroup. Then $H$ is verbally closed in $G$ if and only if $H$ is a retract of~$G$.
\end{thm}

\begin{cor}\label{1_corollary_hyp} {\rm{(\cite[Corollary 2.6]{Bog_1})}
(Solution to Problem 5.2 in~\cite{MR})}\\
Let $G$ be a hyperbolic group.
Suppose that $H$ is a subgroup of $G$ that does not have nontrivial finite normal subgroups
and is not isomorphic to $D_{\infty}$.
Then $H$ is verbally closed in $G$ if and only if $H$ is a retract of $G$.
\end{cor}

The structure of this paper is the following.
Section 2 contains definitions and known results.
In Section~3 we define a composition of two 4-gons.
Section~4 contains a lemma, which is used in Section 5.
In Section 5 we give a proof of a weak version of Theorem~\ref{acylindric},
i.e. a version, where $f(r)$ is, possibly, a non-linear function.
In Section 6 we give a proof of Theorem~\ref{acylindric}.

\medskip

All actions are assumed to be isometric in this paper.
All generating sets considered in this paper are assumed to be symmetric, i.e. closed under taking inverse elements.



\section{Preliminaries}

\subsection{General notation}
Let $G$ be a group generated by a subset $X$. For $g\in G$ let $|g|_X$ be the length of a shortest word in $X$ representing $g$. The corresponding metric
on $G$ is denoted by $d_X$ (or by $d$ if $X$ is clear from the context); thus $ d_X(a,b)=|a^{-1}b|_X$. The right Cayley graph of $G$ with respect to $X$ is denoted by $\Gamma(G,X)$.
By a path $p$ in the Cayley graph we mean a combinatorial path consisting of edges; the initial and the terminal vertices of $p$ are denoted by $p_{-}$ and $p_{+}$, respectively.
The path inverse to $p$ is denoted by $\overline{p}$.
The label of $p$ is denoted by ${\bold {Lab}}(p)$;
we stress that the label is a formal word in the alphabet $X$.
The length of $p$ is denoted by $\ell(p)$.
Given two vertices $a,b$ in $\Gamma(G,X)$, we denote by $[a,b]$ any geodesic path $p$ in $\Gamma(G,X)$ with $p_{-}=a$ and $p_{+}=b$.

Recall that a path $p$ in $\Gamma(G,X)$ is called ($\varkappa,\varepsilon)$-{\it quasi-geodesic} for some $\varkappa\geqslant 1$,
$\varepsilon\geqslant 0$, if $$d(q_{-},q_{+})\geqslant \frac{1}{\varkappa}\ell(q)-\varepsilon$$ for any finite subpath $q$ of $p$.


\subsection{Hyperbolic spaces}

Recall that a metric space $\frak{X}$ is called {\it $\delta$-hyperbolic} if it is geodesic and every geodesic triangle in $\frak{X}$ is $\delta$-slim, i.e. each of the sides of this triangle is contained in the $\delta$-neighborhood of the union of the other two sides (see~\cite[Chapter III.H, Definition~1.1]{BH}).
The following lemma follows straightforward from this definition.

\begin{lem}\label{quadrangle}
Each side of a geodesic quadrangle in a $\delta$-hyperbolic space lies in $2\delta$-neighborhood of the
union of the other three sides.
\end{lem}

We will often use the constant $\mu$ defined in the following lemma.


\begin{lem}\label{close} {\rm (see~\cite[Chapter III.H, Theorem 1.7]{BH})} For all $\delta\geqslant 0$, $\varkappa\geqslant 1$, $\epsilon\geqslant 0$, there exists a constant
$\mu=\mu(\delta,\varkappa,\epsilon)$ with the following property:

If $\frak{X}$ is a $\delta$-hyperbolic space, $p$ is a $(\varkappa,\epsilon)$-quasi-geodesic in $\frak{X}$,
and $[x,y]$ is a geodesic segment joining the endpoints of $p$, then the Hausdorff distance between $[x,y]$
and the image of $p$ is less than $\mu$.
\end{lem}

\begin{cor}\label{close_qg} 
Let $\frak{X}$ be a $\delta$-hyperbolic space, let $p$ and $q$ be $(\varkappa,\epsilon)$-quasi-geodesics in $\frak{X}$ with $\max\{d(p_{-},q_{-}),d(p_{+},q_{+})\}\leqslant r$, and let $\mu=\mu(\delta,\varkappa,\epsilon)$ be the constant from Lemma~\ref{close}. Then the following holds.

\begin{enumerate}
\item[{\rm (1)}] The Hausdorff distance between the images of $p$ and $q$ is less than $r+2(\delta+\mu)$.

\item[{\rm (2)}] If $A\in p$ is a point with
$d(A,\{p_{-},p_{+}\})>r+2\delta+\mu$,
then $A$ lies in the $2(\delta+\mu)$-neighborhood of $q$.
\end{enumerate}
\end{cor}

\medskip
{\it Proof.}
(1) By Lemma~\ref{quadrangle}, the Hausdorff distance between the geodesics $[p_{-},p_{+}]$ and
$[q_{-},q_{+}]$ is at most $(r+2\delta)$.
By Lemma~\ref{close}, the Hausdorff distance between $p$ and $[p_{-},p_{+}]$ (respectively, between $q$ and $[q_{-},q_{+}]$) is less than $\mu$. Adding up, we obtain statement (1).

(2) Let $Q$ be a geodesic quadrangle with vertices $q_{-},p_{-},p_{+},q_{+}$.
The side of $Q$ with endpoints $u$ and $v$ is denoted by $[u,v]$.
Since $p$ lies in the $\mu$-neighborhood of $[p_{-},p_{+}]$, there exists a point $B\in [p_{-},p_{+}]$ such that $d(A,B)\leqslant \mu$. Then $d(B,p_{-})\geqslant d(A,p_{-})-d(A,B)>(r+2\delta+\mu)-\mu=r+2\delta$.
Hence $$d(B,[q_{-},p_{-}])\geqslant d(B,p_{-})-d(p_{-},q_{-})> (r+2\delta)-r=2\delta.$$
Analogously $d(B,[q_{+},p_{+}])> 2\delta$. By Lemma~\ref{quadrangle} applied to $Q$, we have $d(B,[q_{-},q_{+}])\leqslant 2\delta$.
Then $d(A,[q_{-},q_{+}])\leqslant \mu+2\delta$.
Since $[q_{-},q_{+}]$ lies in the $\mu$-neighborhood of $q$, we have $d(A,q)\leqslant 2\mu+2\delta$.
\hfill $\Box$










\begin{lem}\label{quasi-geod} {\rm (see~\cite[Chapter III.H, Theorem 1.13]{BH})}
Let $\frak{X}$ be a $\delta$-hyperbolic geodesic space and let $p:[a,b]\rightarrow \frak{X}$
be a $(8\delta+1)$-local geodesic. Then the following holds:
\begin{enumerate}
\item[(1)] ${\text{\rm im}}(p)$ is contained in the $2\delta$-neighbourhood of any geodesic segment $[p(a),p(b)]$ connecting its endpoints,

\medskip

\item[(2)] $[p(a),p(b)]$ is contained in the $3\delta$-neighborhood of ${\text{\rm im}}(p)$, and

\medskip

\item[(3)] $p$ is a $(3,2\delta)$-quasi-geodesic.
\end{enumerate}
\end{lem}



\medskip

The following definition will help us to shorten the forthcoming proofs.

\begin{defn}
{\rm
For $a,b,c\in \mathbb{R}$, we write $a\approx_{c}b$ if $|a-b|\leqslant c$.
Note that $a\approx_{c}b$ and $b\approx_{c_1}d$ imply $a\approx_{c+c_1}d$.

}
\end{defn}




\begin{lem}\label{approx}
Let $p$ be a $(\varkappa,\epsilon)$-quasi-geodesic in a $\delta$-hyperbolic metric space $\frak{X}$ and let $\mu=\mu(\delta,\varkappa,\epsilon)$ be the constant from Lemma~\ref{close}.
For every three points $A,B,C\in p$ such that $B$ lies on a subpath of $p$ with endpoints $A$ and $C$
we have that $d(A,C)\approx_{2\mu} d(A,B)+d(B,C)$.
\end{lem}

\medskip

{\it Proof.} By Lemma~\ref{close}, there exists a point $B'\in [A,C]$ such that $d(B,B')\leqslant \mu$.
Then the statement follows from $d(A,C)=d(A,B')+d(B',C)$ and $d(A,B')\approx _\mu d(A,B)$, and $d(B',C)\approx _\mu d(B,C)$. \hfill $\Box$

\subsection{Definitions of acylindrically hyperbolic groups}

First we recall a definition of an acylindrical action of a group on a metric space.

\begin{defn} {\rm (see~\cite{Bowditch} and Introduction in~\cite{Osin_2})
An action of a group $G$ on a metric space $\frak{X}$ is called
{\it acylindrical}
if for every $\varepsilon\geqslant 0$ there exist $R,N>0$ such that for every two points $x,y$ with $d(x,y)\geqslant R$,
there are at most $N$ elements $g\in G$ satisfying
$$
d(x,gx)\leqslant \varepsilon\hspace*{2mm}{\text{\rm and}}\hspace*{2mm} d(y,gy)\leqslant \varepsilon.
$$
}
\end{defn}


For Cayley graphs, the acylindricity condition can be rewritten as follows:


\begin{defn}\label{A}
{\rm Let $G$ be a group and $X$ be a generating set of $G$. The natural action of $G$
on the Cayley graph $\Gamma(G,X)$ is called {\it acylindrical}
if for every $\varepsilon\geqslant 0$ there exist $R,N>0$ such that for any $g\in G$ of length $|g|_X\geqslant R$
we have
$$
\bigl|\{f\in G\,|\, |f|_X\leqslant \varepsilon,\hspace*{2mm} |g^{-1}fg|_X\leqslant \varepsilon \}\bigr|\leqslant N.
$$
}
\end{defn}

\medskip

Recall that an action of a group $G$ on a hyperbolic space $S$ is called {\it elementary} if the limit set
of $G$ on the Gromov boundary $\partial S$ contains at most 2 points.

\begin{defn}\label{definition_acylindrically_hyperbolic} {\rm (see~\cite[Definition 1.3]{Osin_2})
A group $G$ is called {\it acylindrically hyperbolic} if it satisfies one of the following equivalent
conditions:

\begin{enumerate}
\item[(${\rm AH}_1$)] There exists a generating set $X$ of $G$ such that the corresponding Cayley graph $\Gamma(G,X)$
is hyperbolic, $|\partial \Gamma (G,X)|>2$, and the natural action of $G$ on $\Gamma(G,X)$ is acylindrical.

\medskip

\item[(${\rm AH}_2$)] $G$ admits a non-elementary acylindrical action on a hyperbolic space.
\end{enumerate}
}
\end{defn}

In the case (${\rm AH}_1$), we also write that $G$ is {\it acylindrically hyperbolic with respect to $X$}.

\subsection{Elliptic and loxodromic elements in acylindrically hyperbolic groups}

Given a group $G$ acting on a metric space $\frak{X}$ by isometries, an element $g\in G $ is called {\it elliptic} if some (equivalently, any) orbit of $g$ is bounded, and
{\it loxodromic} if the map $\mathbb{Z}\rightarrow \frak{X}$ defined by $n\mapsto g^n(x)$  is a quasi-isometric embedding for some (equivalently, any) $x\in \frak{X}$. Even in the case of groups acting on hyperbolic spaces, there may be other types of actions (see~\cite[Section 8.2]{Gromov}  and~\cite[Section 3]{Osin_2}).

Let $X$ be a generating set of $G$. An element $g\in G$ is called {\it loxodromic (elliptic) with respect to $X$} if $g$ is loxodromic (elliptic) for the left action of $G$ on the Cayley graph $\Gamma(G,X)$.
In other words, $g\in G$ is loxodromic with respect to $X$ if and only if
there exist $\eta>0$ and $\nu\geqslant 0$ such that $|g^n|_X\geqslant \eta n-\nu $ for all $n\in \mathbb{N}$.

Equivalently, an element $g\in G$ is loxodromic with respect to $X$ if the path $L(1,g)$ in $\Gamma(G,X)$
is a $(\varkappa,\epsilon)$-quasi-geodesic for some $\varkappa\geqslant 1$ and $\epsilon\geqslant 0$.




Recall that the {\it stable norm} of an element $g\in G$ with respect to $X$ is defined~as $$||g||_X=\underset{n\rightarrow \infty}{\lim}\frac{|g^n|_X}{n},$$
see~\cite{CDP}. This limit exists since the function $f:\mathbb{N}\rightarrow \mathbb{N}\cup \{0\}$, $n\mapsto |g^n|_X$, is subadditive.
The following lemma is quite trivial.

\begin{lem}\label{inf} The stable norm satisfies the following properties.

\begin{enumerate}
\item[{\rm (1)}] $\displaystyle{|g|_X\geqslant ||g||_X=\underset{n\in \mathbb{N}}{\inf}\frac{|g^n|_X}{n}}$.

\medskip

\item[{\rm (2)}] $||y^{-1}gy||_X=||g||_X$ for any $y\in G$.

\medskip

\item[{\rm (3)}] $||g^k||_X=|k|\cdot ||g||_X$ for any $k\in \mathbb{Z}$.
\end{enumerate}
\end{lem}

\begin{lem}\label{Bowditch} {\rm (see~\cite[Lemma 2.2]{Bowditch})}
Let $X$ be a generating set of a group $G$.\break If the Cayley graph $\Gamma(G,X)$ is hyperbolic and $G$ acts acylindrically on $\Gamma(G,X)$, then each element of $G$ is either elliptic or loxodromic with respect to $X$.
Moreover, there exists $\tau>0$ such that
the stable norm of any loxodromic element of $G$ is at least $\tau$.
\end{lem}

\begin{cor}\label{qg} Let $G$ be a group and let $X$ be a generating set of $G$. Suppose that the Cayley graph $\Gamma(G,X)$ is $\delta$-hyperbolic and that $G$ acts acylindrically on $\Gamma(G,X)$.
Then there exist $\varkappa_0\geqslant 1$ and $\epsilon_0\geqslant 0$ such that the following holds:
If an element $g\in G$ is loxodromic and shortest in its conjugacy class, then the quasi-geodesic $L(1,g)$
associated with $g$ is a $(\varkappa_0,\epsilon_0)$-quasi-geodesic.
\end{cor}


{\it Proof.} We consider two cases.

{\bf Case 1.} Suppose that $|g|_X\geqslant 8\delta+1$.

Since $g$ is shortest in its conjugacy class, the path $L(1,g)$ is a $(8\delta+1)$-local geodesic. Then, by Lemma~\ref{quasi-geod},
$L(1,g)$ is a $(3,2\delta)$-quasi-geodesic.


{\bf Case 2.} Suppose that $|g|_X<8\delta+1$.

By Lemmas~\ref{inf} (1) and~\ref{Bowditch}, there exists a constant $\tau>0$ which does not depend on $g$
and such that
$|g^n|_X\geqslant n\tau$  for any $n\in \mathbb{N}$. In particular.
$$
8\delta+1>|g|_X\geqslant \tau.\eqno{(2.1)}
$$
We prove that $L(1,g)$ is a $(\sigma,\varepsilon)$-quasi-geodesic for $\sigma=\frac{8\delta+1}{\tau}$
and $\varepsilon=4(8\delta+1)$.
Let $p$ be a subpath of $L(1,g)$. If $p$ does not contain phase vertices of $L(1,g)$,
then $p$ is a geodesic path, and hence a $(1,0)$-quasi-geodesic.
Suppose that $p$ contains at least one phase vertex of $L(1,g)$. Let $q$ be the maximal subpath of $p$ such that
$q_{-}=g^i$ and $q_{+}=g^{i+n}$ for some $i\in \mathbb{Z}$ and $n\in \mathbb{N}\cup \{0\}$. Then
$$
\begin{array}{ll}
d_X(p_{-},p_{+}) & \geqslant d_X(q_{-},q_{+})-2|g|_X\vspace*{2mm}\\
& =|g^n|_X-2|g|_X \vspace*{2mm}\\
& \geqslant n\tau-2|g|_X\vspace*{2mm}\\
& > \frac{1}{\sigma}n|g|_X - 2|g|_X\vspace*{2mm}\\
& \geqslant \frac{1}{\sigma}(\ell(p)-2|g|_X)-2|g|_X> \frac{1}{\sigma}\ell(p)- \varepsilon,
\end{array}
$$
since
$$
\begin{array}{l}
(\frac{1}{\sigma}+1)2|g|_X=(\frac{\tau}{8\delta+1}+1)2|g|_X\overset{(2.1)}{<}4(8\delta+1)=\varepsilon.
\end{array}
$$
Thus, $L(1,g)$ is a $(\sigma,\varepsilon)$-quasi-geodesic in this case.

In both cases $L(1,g)$ is a $(\varkappa_0,\epsilon_0)$-quasi-geodesic for $\varkappa_0=\max\{3,\sigma\}$
and $\epsilon_0=4(8\delta+1)$.
\hfill $\Box$

\medskip

Recall that any loxodromic element $g$ in an acylindrically hyperbolic group $G$ is contained in a
unique maximal virtually cyclic subgroup $E_G(g)$ of $G$, see~\cite[Lemma 6.5]{DOG}. This subgroup is called the {\it elementary subgroup associated with $g$}; it can be described as follows (see ~\cite[Corollary 6.6]{DOG}):
$$E_G(g)=\{f\in G\,|\, \exists  n\in \mathbb{N}:  f^{-1}g^nf=g^{\pm n}\}.$$
The centralizer of $g$ in $G$, denoted by $C_G(g)$, is contained in $E_G(g)$.





\section{Composition of 4-gons}

With the term {\it $n$-gon} in $\Gamma(G,X)$, we mean an ordered sequence of paths $p_1,\dots,p_n$ in
$\Gamma(G,X)$ such that $(p_i)_{+}=(p_{i+1})_{-}$ for $i=1,\dots ,n-1$, and $(p_n)_{+}=(p_1)_{-}$. The paths $p_i$ are called {\it sides} of this $n$-gon. Thus, the sides are
not necessarily geodesics, and they can intersect each other.

Given a 4-gon $P=p_1p_2p_3p_4$, we say that $p_1,p_2,\overline{p_3}$ and $\overline{p_4}$ are the left,
the top, the right, and the bottom sides of $P$, respectively. By ${\bf L}(P)$, ${\bf T}(P)$, ${\bf R}(P)$,
${\bf B}(P)$ denote the elements of $G$ corresponding to the labels of the left, top, right, and bottom
sides of $P$.

Let $P=p_1p_2p_3p_4$ and $Q=q_1q_2q_3q_4$ be two 4-gons in $\Gamma(G,X)$ such that
their top sides have the same label. Then there exists a unique element $g\in G$ such that $g(q_2)=p_2$.
Consider the 4-gon $g(Q)=r_1r_2r_3r_4$, where $r_i=g(q_i)$. Since the top sides of $P$ and $g(Q)$ coincide,
we can glue these polygons along their top sides and obtain a new 4-gon $S=(p_1\overline{r_1})\overline{r_4}(\overline{r_3}p_3)p_4$ in $\Gamma(G,X)$. In particular, $p_1\overline{r_1}$ is the left side of $S$.
We say that $S$ is the {\it composition} of $P$ and $Q$ and denote $S$ by $P\circ Q$, see Fig. 3.

\medskip

\vspace*{-32mm}
\hspace*{-10mm}
\includegraphics[scale=0.7]{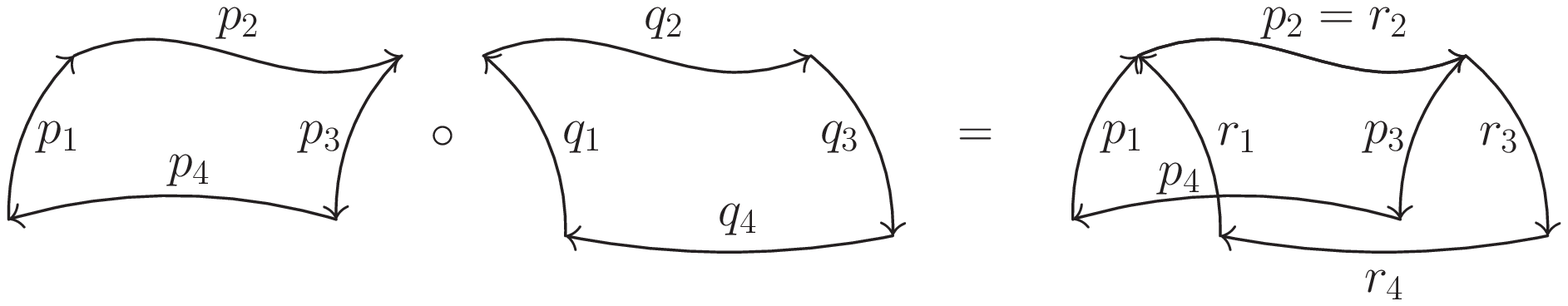}

\vspace*{-14.6cm}

\begin{center}
Fig. 3. Composition of two 4-gons.
\end{center}

\medskip

Let $\psi:Q\rightarrow P\circ Q$ be the label preserving map which sends $v$ to $gv$ for any vertex $v$ of $Q$.
We call $\psi$ the {\it translation map} for the pair $(P,Q)$.

\newpage

\section{Auxiliary lemma}

Let $b\in G$. We say that a quasi-geodesic path $\mathcal{P}$ in $\Gamma(G,X)$ is {\it $b$-periodic}
if $\mathcal{P}=L(x,b)$ for some $x\in G$.

\begin{lem}\label{conjugacy_2}
Let $G$ be a group and $X$ a generating set of $G$. Suppose that the Cayley graph $\Gamma(G,X)$ is hyperbolic and that $G$ acts acylindrically on $\Gamma(G,X)$.
Then there is a function $K:\mathbb{N}\cup\{0\}\rightarrow \mathbb{N}$ such that the following holds.

Let $b\in G$ be an element that is loxodromic (with respect to $X$) and shortest in its conjugacy class.
Let $\mathcal{P}$ and $\mathcal{Q}$ be two $b$-periodic quasi-geodesic paths in $\Gamma(G,X)$ and $r$ a non-negative integer.
Suppose that $p$ and $q$ are finite subpaths of  $\mathcal{P}$ and $\mathcal{Q}$ such that
$$\max\{d(p_{-}, q_{-}),\,d(p_{+}, q_{+})\}\leqslant r,$$
and $p$ and $q$ contain at least $K(r)$ $b$-periods of $\mathcal{P}$ and $\mathcal{Q}$, respectively.
Then for any phase vertices $B\in \mathcal{P}$ and $A\in \mathcal{Q}$,
the element $A^{-1}B$ centralizes a nontrivial power of $b$.
\end{lem}

{\it Proof.} First we define some constants which we use in the proof.

\begin{enumerate}
\item[$\bullet$] We fix $\delta\geqslant 0$ such that the Cayley graph $\Gamma(G,X)$ is $\delta$-hyperbolic.

\item[$\bullet$] By Corollary~\ref{qg}, there exist universal constants $\varkappa_0\geqslant 1$ and $\epsilon_0\geqslant 0$ such that the paths $\mathcal{P}$ and $\mathcal{Q}$ are $(\varkappa_0,\epsilon_0)$-quasi-geodesic.

\item[$\bullet$] Let $\mu=\mu(\delta,\varkappa_0,\epsilon_0)$ be the constant from Lemma~\ref{close}.

\item[$\bullet$] We set $\varepsilon:=6r+24\mu+8\delta$.
Let $R=R(\varepsilon)$ and $N=N(\varepsilon)$ be constants from Definition~\ref{A}.
We may assume that $N$ is a positive integer.
We also set $S:=\lceil\varkappa_0(R+\epsilon_0)\rceil$ and define
$$
K(r)=S+N+1.
$$
\end{enumerate}


We show that $K(r)$ is the desired number. Thus, suppose that $p$ and $q$ contain at least $K(r)$
$b$-periods of $\mathcal{P}$ and $\mathcal{Q}$, respectively.
Let $B_0$ be the leftmost phase vertex of $\mathcal{P}$ lying in $p$ and let $A_0$ be the leftmost phase vertex of $\mathcal{Q}$ lying in $q$  (see Fig.~4).
By assumption, the vertices $B_j:=B_0b^j$ and  $A_j:=A_0b^j$, $0\leqslant j\leqslant K(r)$,
are phase vertices lying on $p$ and $q$, respectively.
Let $u$ be the segment of $p$ from $p_{-}$ to $B_0$, and let $v$ be the segment of $q$ from $q_{-}$ to $A_0$. Without loss of generality we assume that $\ell(u)\geqslant \ell(v)$. Note that $\ell(u)<|b|_X$.

For $1 \leqslant i\leqslant K(r)$, let $C_i$ be the vertex on $[B_{i-1},B_i]$ such that $d(C_i,B_i)=\ell(u)-\ell(v)$.
This definition implies the following claim.

\medskip

{\bf Claim 1.} The element $C_i^{-1}B_i$ does not depend on $i$.

\medskip
\vspace*{-35mm}
\hspace*{-10mm}
\includegraphics[scale=0.8]{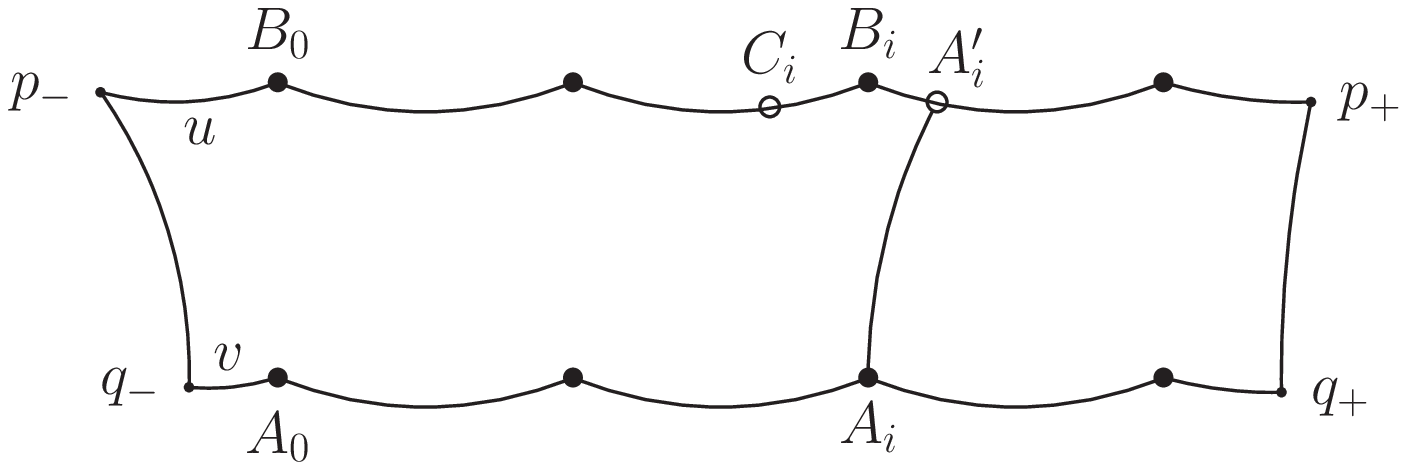}

\vspace*{-16cm}

\begin{center}
Fig. 4. Illustration to Claims 1 and 2.\\
\end{center}

\medskip

{\bf Claim 2.}  $d(A_i,C_i)\leqslant \varepsilon/2.$

\medskip

{\it Proof.} By Lemma~\ref{approx}, we have
$$
d(p_{-},B_i)\approx_{2\mu} d(p_{-},B_0)+d(B_0,B_i)=d(B_0,B_i)+\ell(u)
$$
and
$$
d(p_{-},B_i)\approx_{2\mu} d(p_{-},C_i)+d(C_i,B_i)=d(p_{-},C_i)+\ell(u)-\ell(v).
$$
Therefore
$$
d(p_{-},C_i)\approx_{4\mu} d(B_0,B_i)+\ell(v).
$$
On the other hand
$$
d(q_{-},A_i)\approx_{2\mu} d(q_{-},A_0)+d(A_0,A_i)=d(A_0,A_i)+\ell(v)=d(B_0,B_i)+\ell(v).
$$
Therefore
$$
d(p_{-},C_i)\approx_{6\mu} d(q_{-},A_i).\eqno{(4.1)}
$$

\medskip

Since $\max\{d(p_{-},q_{-}), d(p_{+},q_{+})\}\leqslant r$, there exists a vertex $A_i'\in p$ such that
$$
d(A_i,A_i')\leqslant r+2\mu+2\delta
$$
(see Corollary~\ref{close_qg}). Considering any geodesic 4-gon with vertices $p_{-}$, $A_i'$, $A_i$, $q_{-}$, we obtain
$$|d(p_{-},A_i')- d(q_{-},A_i)|\leqslant d(p_{-},q_{-})+d(A_i,A_i')\leqslant r+(r+2\mu+2\delta),$$
i.e.
$$d(p_{-},A_i')\approx_{2(r+\mu+\delta)}d(q_{-},A_i).\eqno{(4.2)}$$
It follows from (4.1) and (4.2) that
$$
|d(p_{-},C_i)-d(p_{-},A_i')|\leqslant 2r+8\mu+2\delta.\eqno{(4.3)}
$$
On the other hand, by Lemma~\ref{approx}, we have
$$
d(p_{-},C_i)-d(p_{-},A_i')\approx_{2\mu} d(C_i,A_i')
$$
if $C_i$ lies on $p$ to the right of $A_i'$ and
$$
d(p_{-},A_i')-d(p_{-},C_i)\approx_{2\mu} d(C_i,A_i')
$$
if $A_i'$ lies on $p$ to the right of $C_i$.
In each of these cases, we deduce from (4.3) that
$$
d(C_i,A_i')\leqslant 2r+10\mu+2\delta.
$$

Finally,
$
d(C_i,A_i)\leqslant d(C_i,A_i')+d(A_i',A_i)\leqslant 3r+12\mu+4\delta=\varepsilon/2.
$
\hfill $\Box$

\bigskip

Recall that $S:=\lceil\varkappa_0(R+\epsilon_0)\rceil$.
For $i=1,\dots, K(r)-S$, we consider the 4-gons $\mathcal{P}_i=p_{1,i}p_{2,i}p_{3,i}p_{4,i}$,
where $p_{1,i}$ is a geodesic from $A_i$ to $B_i$,
$p_{2,i}$ is the subpath of $p$ from $B_i$ to $B_{i+S}$,
$p_{3,i}$ is a geodesic from $B_{i+S}$ to $A_{i+S}$,
and $p_{4,i}$ is the subpath of $\bar{q}$ from $A_{i+S}$ to $A_i$.

Since the top sides of $\mathcal{P}_1$ and $\mathcal{P}_i$ have the same labels,
we can build the composition $\mathcal{P}_1\circ \mathcal{P}_i$.

\medskip

{\bf Claim 3.} For $1\leqslant i\leqslant K(r)-S$ the following hold.
\begin{enumerate}
\item[(a)] 
${\bf T}(\mathcal{P}_1\circ \mathcal{P}_i)=b^S$ and ${\bf B}(\mathcal{P}_1\circ \mathcal{P}_i)=b^S$.

\medskip

\item[(b)]
$|{\bf L}(\mathcal{P}_1\circ \mathcal{P}_i)|_X\leqslant \varepsilon$ and $|{\bf R}(\mathcal{P}_1\circ \mathcal{P}_i)|_X\leqslant \varepsilon$.

\end{enumerate}

\medskip

{\it Proof.}
Statement (a) follows from definition. We prove (b).
Using Claim 1, we deduce
$$
\begin{array}{ll}
{\text{\bf L}}(\mathcal{P}_1\circ \mathcal{P}_i) & = {\text {\bf L}}(\mathcal{P}_1)({\bf L}(\mathcal{P}_i))^{-1}\vspace*{2mm}\\
 & = (A_1^{-1}B_1)(A_i^{-1}B_i)^{-1}\vspace*{2mm}\\
 & = (A_1^{-1}C_1)(C_1^{-1}B_1)(B_i^{-1}C_i)(C_i^{-1}A_i)\vspace*{2mm}\\
 & = (A_1^{-1}C_1)(C_i^{-1}A_i).
\end{array}
$$
Hence, by Claim 2, we have
$$
|{\text{\bf L}}(\mathcal{P}_1\circ \mathcal{P}_i)|_X\leqslant d(A_1,C_1)+d(A_i,C_i)\leqslant \varepsilon.
$$
Analogously, we have
$
|{\text{\bf R}}(\mathcal{P}_1\circ \mathcal{P}_i)|_X\leqslant \varepsilon.
$
\hfill $\Box$

\bigskip

Since $G$ is acylindrically hyperbolic with respect to $X$ and
$$
|b^S|_X\geqslant \frac{1}{\varkappa_0}S|b|_X-\epsilon_0\geqslant \frac{1}{\varkappa_0}S-\epsilon_0\geqslant  R,
$$
it follows from Claim 3 that the number of elements
${\text{\bf L}}(\mathcal{P}_1\circ \mathcal{P}_i)$, where $i$ runs from 1 to $K(r)-S=N+1$,
is at most $N$.
Therefore there are different $k,l$ such that $${\text{\bf L}}(\mathcal{P}_k)={\text{\bf L}}(\mathcal{P}_l).\eqno{(4.4)}$$
We denote this element by $z$.
We have
$$
z=A_k^{-1}B_k\overset{(4.4)}{=}A_l^{-1}B_l=(A_l^{-1}A_k)(A_k^{-1}B_k)(B_k^{-1}B_l)
=b^{k-l}zb^{l-k},$$
hence $z$ centralizes $b^{k-l}$.

Finally, suppose that $A$ and $B$ are arbitrary phase vertices of $\mathcal{Q}$ and $\mathcal{P}$, respectively.
Then $A_k^{-1}A$ and $B_k^{-1}B$ lie in $\langle b\rangle$. It follows that
$A^{-1}B$ also centralizes $b^{k-l}$.
\hfill $\Box$

\section{A weak version of Theorem~\ref{acylindric}}

We show that Theorem~\ref{acylindric} is valid for the following (possibly nonlinear) function $F(r)$ instead of $f(r)$: $$F(r)=\varkappa_0\bigl(K(2r)+(\epsilon_0+2r+2)\bigr)+1,\eqno{(5.1)}$$ where $\varkappa_0$ and $\epsilon_0$ are the numbers from Corollary~\ref{qg} and
$K$ is the function from Lemma~\ref{conjugacy_2}. Increasing $\varkappa_0$ and $\epsilon_0$, we may assume that
these numbers are integers. Thus, $F(r)\in \mathbb{N}$ for $r\in \mathbb{N}\cup \{0\}$.

\begin{lem}\label{weak_acylindric}
Let $G$ be a group and $X$ a generating set of $G$. Suppose that the Cayley graph $\Gamma(G,X)$
is hyperbolic and that $G$ acts acylindrically on $\Gamma(G,X)$.
Then the following holds.

Let $a,b\in G$ be two loxodromic elements which are shortest in their conjugacy classes
and such that $|a|_X\geqslant |b|_X$.
Let $x,y\in G$ be arbitrary elements and $r$ an arbitrary non-negative integer.

If a subpath $p\subset L(x,a)$ contains at least $F(r)$ $a$-periods and lies in the $r$-neighborhood of $L(y,b)$,
then there exist $s,t\neq 0$ such that $(x^{-1}y)b^s(y^{-1}x)=a^t$.
\end{lem}


\medskip

{\it Proof.}
Let $A_0,\dots ,A_{F(r)},\dots $ be the consequent phase vertices of $L(x,a)$ lying on $p$.
By assumption, for each $0\leqslant i\leqslant F(r)$, there exists a vertex $C_i\in L(y,b)$
such that $d(A_i,C_i)\leqslant r$.
Let $m_i$ be the minimal integer such that $C_i$ lies on the subpath $[yb^{m_i-1},yb^{m_i}]$ of $L(y,b)$. We denote $D_i=yb^{m_i}$.
Thus, $D_i$ is the nearest phase vertex of $L(y,b)$ which lies to the right of $C_i$ (see Fig.~5).

\medskip

\vspace*{-30mm}
\hspace*{-12mm}
\includegraphics[scale=0.7]{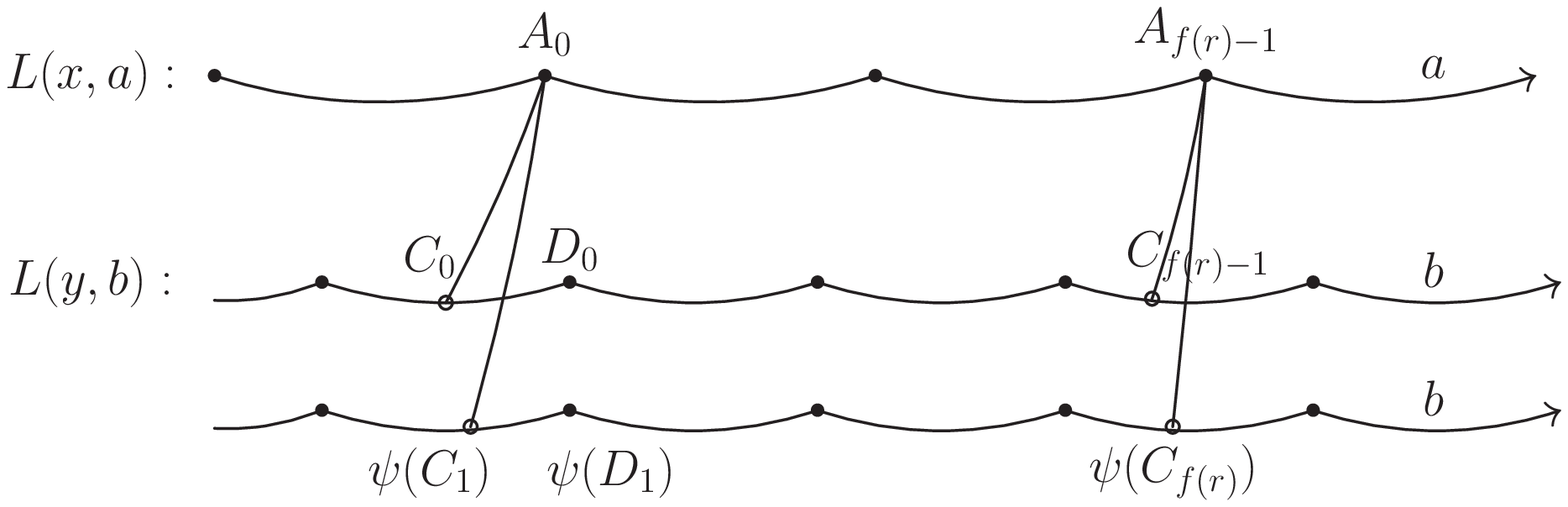}

\vspace*{-125mm}

\begin{center}
Fig. 5. Illustration to the proof of Lemma~\ref{weak_acylindric}.
\end{center}


For $i=0,1$, we consider the 4-gon $\mathcal{Q}_i=q_{1,i}q_{2,i}q_{3,i}q_{4,i}$,
where $q_{1,i}$ is a geodesic from $C_i$ to $A_i$,
$q_{2,i}$ is the subpath of $L(x,a)$ from $A_i$ to $A_{i+F(r)-1}$,
$q_{3,i}$ is a geodesic from $A_{i+F(r)-1}$ to $C_{i+F(r)-1}$,
and $q_{4,i}$ is the subpath of $L(y,b)$ from $C_{i+F(r)-1}$ to $C_i$.

Since the top sides of $\mathcal{Q}_0$ and $\mathcal{Q}_1$ have the same label (representing the element $a^{F(r)-1}$), the composition $\mathcal{Q}_0\circ \mathcal{Q}_1$ is defined. Let $\psi:\mathcal{Q}_1\rightarrow \mathcal{Q}_0\circ \mathcal{Q}_1$
be the translation map (see Section 3) for the pair $(\mathcal{Q}_0,\mathcal{Q}_1)$.
The map $\psi$ is determined by the element $g\in G$ such that $A_0=gA_1$. We have $A_0=xa^i$ and $A_1=xa^{i+1}$ for some $i$. Then  $g=A_0A_1^{-1}=xa^{-1}x^{-1}$.

\medskip

{\bf Claim.}
\begin{enumerate}

\item[1)] The bottom and the top sides of $\mathcal{Q}_0\circ \mathcal{Q}_1$ are parts of $b$-periodic quasi-geodesics $L(y,b)$ and $(xa^{-1}x^{-1})L(y,b)$, respectively.

\item[2)] The length of the left and of the right sides of the 4-gon $\mathcal{Q}_0\circ \mathcal{Q}_1$ is at most $2r$.

\item[3)] The bottom and the top sides of $\mathcal{Q}_0\circ \mathcal{Q}_1$ contain at least $K(2r)$\break $b$-periods.
\end{enumerate}

{\it Proof.} Statement 1) follows from definition.

2) The length of the left side of $\mathcal{Q}_0\circ \mathcal{Q}_1$ equals the sum of lengths of the left sides of $\mathcal{Q}_0$ and $\mathcal{Q}_1$, i.e. is equal to $d(C_0,A_0)+d(C_1,A_1)\leqslant 2r$.

3) Let $N$ be the number of $b$-periods on the bottom side $\overline{q_{4,0}}$. We have
$$
\begin{array}{ll}
(N+2)|b|_X\geqslant \ell(q_{4,0}) & \geqslant d(C_0, C_{F(r)-1})\vspace*{2mm}\\
& \geqslant d(A_0, A_{F(r)-1})-d(C_0,A_0)-d(C_{F(r)-1}, A_{F(r)-1})\vspace*{2mm}\\
& \geqslant d(A_0, A_{F(r)-1})-2r\vspace*{2mm}\\
& \geqslant \frac{1}{\varkappa_0} (F(r)-1)|a|_X-\epsilon_0-2r. \vspace*{2mm}\\
\end{array}
$$

This implies
$$
N\geqslant \frac{1}{\varkappa_0}(F(r)-1)\frac{|a|_X}{|b|_X}-(\epsilon_0+2r+2)\geqslant
\frac{1}{\varkappa_0}(F(r)-1)-(\epsilon_0+2r+2)= K(2r).
$$

Analogously, the number of $b$-periods on the top side $g(\overline{q_{4,1}})$ of
$\mathcal{Q}_0\circ \mathcal{Q}_1$ is at least $K(2r)$.\hfill $\Box$

\medskip

This claim and Lemma~\ref{conjugacy_2} imply
$$
D_0^{-1}(\psi(D_1))\in C_G(b^n).\eqno{(5.2)}
$$
for some $n\in \mathbb{N}$. We have
$$
\begin{array}{ll}
D_0^{-1}(\psi(D_1)) & = (D_0^{-1}A_0)(A_0^{-1}\psi(D_1))\vspace*{2mm}\\
                    & = (D_0^{-1}A_0)((\psi(A_1))^{-1}\psi(D_1))\vspace*{2mm}\\
                    & = (D_0^{-1}A_0)(A_1^{-1} D_1).\\
\end{array}\eqno{(5.3)}
$$

Recall that $A_0=xa^i$ and $A_1=xa^{i+1}$ for some $i\in \mathbb{Z}$, and that $D_0=yb^{j_0}$, $D_1=yb^{j_1}$ for some $j_0,j_1\in \mathbb{Z}$.
Then, with the help of (5.2) and (5.3), we obtain
$$
(D_0^{-1}A_0)(A_1^{-1}D_1)=b^{-j_0}(y^{-1}x)a^{-1}(x^{-1}y)b^{j_1}\in C_G(b^n).
$$
This implies $(y^{-1}x)a(x^{-1}y)\in C_G(b^n)$.
Since $C_G(b^n)$ lies in the virtually cyclic group $E_G(b)$,
we have $|C_G(b^n):\langle b\rangle|<\infty$.
Therefore there exist $s,t\in \mathbb{Z}\setminus \{0\}$
such that $(y^{-1}x)a^t(x^{-1}y)=b^s$.
\hfill $\Box$

\section{Proof of Theorem~\ref{acylindric}}

We set $\tau={\bf{inj}}(G,X)$ and $\mu=\mu(\delta,\varkappa_0,\epsilon_0)$, where $\mu$ is the function from Lemma~\ref{close} and $\varkappa_0$ and $\epsilon_0$ are the numbers from Corollary~\ref{qg}.
Increasing $\mu$, we may assume that $2\delta+2\mu\in \mathbb{N}\cup \{0\}$.
We show that in Theorem~\ref{acylindric} one can take the constant
$$
C=F(2\delta+2\mu)+\frac{6\mu+4\delta}{\tau}+2,
$$
where $F$ is the function in (5.1).
Let $$
f(r)=\frac{2r}{\tau}+C
$$
and suppose that $p$ is a subpath of $L(x,a)$ which contains at least $f(r)$ $a$-periods
and lies in the $r$-neighborhood of a finite subpath $q$ of $L(y,b)$.
We shall prove that $(x^{-1}y)b^s(y^{-1}x)=a^t$ for some $s,t\neq 0$.

Let $A_0,\dots ,A_n$ be the consequent phase vertices of $L(x,a)$ lying on $p$. By assumption $n\geqslant f(r)$.
Let $p_{k,n-k}$ be the subpath of $p$ from $A_k$ to $A_{n-k}$, where
$$k=\Bigl\lfloor\frac{f(r)-F(2\delta+2\mu)}{2}\Bigr\rfloor=
\Bigl\lfloor\frac{r+3\mu+2\delta}{\tau}+1\Bigr\rfloor
.\eqno{(6.1)}$$
Then Theorem~\ref{acylindric} follows from Lemma~\ref{weak_acylindric} applied to the subpath $p_{k,n-k}\subset L(x,a)$ and the following claim.

\medskip

{\bf Claim.}
The path $p_{k,n-k}\subset L(x,a)$ satisfies the following properties.

1) $p_{k,n-k}$ contains at least $F(2\mu+2\delta)$ $a$-periods.

2) $p_{k,n-k}$ lies in the $(2\mu+2\delta)$-neighborhood of the subpath $q\subset L(y,b)$.

\medskip

{\it Proof.} 1) Since $p$ contains at least $f(r)$ $a$-periods,
the path $p_{k,n-k}$ contains at least $f(r)-2k\geqslant F(2\delta+2\mu)$\, $a$-periods.

2) Let $A\in p_{k,n-k}$ be an arbitrary vertex. Note that $A\in p$. By Corollary~\ref{close_qg}, it suffices to prove that
$d(p_{-},A)>r+2\delta+\mu\hspace*{2mm}{\text{\rm and}}\hspace*{2mm}d(p_{+},A)>r+2\delta+\mu.$
Because of symmetry, we prove only the first inequality. Let $\widetilde{[p_{-},A]}$ be the part of the quasi-geodesic $p$ from $p_{-}$ to $A$.
Since $A_0,A_k\in \widetilde{[p_{-},A]}$, there exist points $A_0',A_k'\in [p_{-},A]$
such that $d(A_0,A_0')\leqslant \mu$ and $d(A_k,A_k')\leqslant \mu$. Then
$$
d(p_{-},A)\geqslant d(A_0',A_k')\geqslant d(A_0,A_k)-2\mu\geqslant
|a^k|_X-2\mu\geqslant \tau k-2\mu\overset{(6.1)}{>}r+2\delta+\mu.
$$
\hfill $\Box$


\end{document}